\documentclass{elsart}

\journal{Journal Mathematical Analysis and Applications}

\usepackage{amsmath}
\usepackage{amsfonts}
\usepackage{amscd}

\theoremstyle{plain} 

\newtheorem{theorem}{Theorem}[section]

\newtheorem{proposition}[theorem]{Proposition}

\newtheorem{remark}[theorem]{Remark}

\begin{document}

\begin{frontmatter}


\title{A note on existence and uniqueness of limit cycles for Li\'enard systems.}

\author[tim]{Timoteo Carletti\corauthref{cor}},
\corauth[cor]{Corresponding author.}
\ead{t.carletti@sns.it}
\author[vil]{Gabriele Villari}
\ead{villari@math.unifi.it}

\address[tim]{Scuola Normale Superiore, piazza dei Cavalieri 7, 56126  Pisa, Italy}
\address[vil]{Dipartimento Matematica U. Dini, viale Morgagni 67/A, 50134 Firenze, Italy}

\begin{abstract}
We consider the Li\'enard equation and we give a sufficient
condition to ensure existence and uniqueness of limit cycles. We
compare our result with some other existing ones and we give some
applications.
\end{abstract}

\end{frontmatter}

\section{Introduction}
\label{sec:intro}

In this paper we consider the Li\'enard equation:
\begin{equation}
\label{eq:lienard1}
\ddot x +f(x) \dot x + g(x)=0 \, ,
\end{equation}
where $f,g:\mathbb{R}\rightarrow \mathbb{R}$, with particular attention
to the existence and uniqueness of limit cycles. This is a classical
problem of non--linear oscillation for second order differential
equations. Different assumptions on $f$ and $g$ and different methods
used to study the problem, gave rise to a large amount of literature
on this topic; for a review of results and methods, reader can
consult~\cite{Staude81}, Chapter IV of the
book~\cite{Zhang92} or~\cite{Ye}. In the following we will give some more references.

We make the following assumptions on $f$ and $g$:
\begin{itemize}
\item[(A)] $f$ is a continuous function and $g$ verifies a locally Lipschitz condition;
\item[(B)] $f(0)<0$, $f(x)>0$ for $|x|>\delta$, for some $\delta >0$, and $xg(x)>0$ for $x\neq 0$.
\end{itemize}
For the Li\'enard equation condition (A) assures existence and uniqueness of the Cauchy
initial value problem, in fact passing to the {\em Li\'enard plane} the second order differential equation is {\em equivalent} to
the following first order system:
\begin{equation}
\label{eq:lienard2}
\begin{cases}
\dot x &= y-F(x) \\
\dot y &= -g(x) \, ,
\end{cases}
\end{equation}
where $F(x)=\int_0^x f(\xi) \, d\xi$. Hence assuming hypothesis (A) the right hand side of~\eqref{eq:lienard2} is Lipschitz
continuous, from which the claim follows. 

Assumption (B) guarantees that the 
origin is the only singular point of the system, which results a repellor, moreover orbits of~\eqref{eq:lienard2} 
turn clockwise around it. Hypothesis on the sign of $f(0)$ can be weakened by asking $xF(x)<0$ for $|x|$ 
small, we nevertheless prefer the former formulation (A) because of the applications we will give in the last 
part of the paper.

Assuming assumptions (A) and (B) on $f$ and $g$, our main result will be the following

\begin{theorem}
\label{thm:main}
Let $G(x)=\int_0^x g(\xi) \, d\xi$ and suppose that $F$ and $G$ verify:
\begin{itemize}
\item[(C)] $F$ has only three real transversal zeros, located at
 $x_0=0$, $x_2<0<x_1$. Assume moreover that $F$ is monotone increasing outside
the interval $[x_2,x_1]$;
\item[(D)] $G(x_1) = G(x_2)$;
\item[(E)] $\limsup_{x\rightarrow +\infty}[G(x)+F(x)]=+\infty$ and 
$\limsup_{x\rightarrow -\infty}[G(x)-F(x)]=+\infty$.
\end{itemize}
Then system~\eqref{eq:lienard2} has a unique periodic orbit in the
 $(x,y)$--plane which is stable.
\end{theorem}
Because of the equivalence of equation~\eqref{eq:lienard1} and system~\eqref{eq:lienard2} also the former has a unique limit cycle
if Theorem applies.

We postpone the proof of the Theorem to the next section, in the
 following we will discuss the role of our hypotheses 
 and compare this result with other existence and uniqueness results
 concerning Lin\'enard systems.

Our result follows from investigating 
the geometry of limit cycles, in particular their (eventual) intersections
with the lines $x=x_1$ and $x=x_2$. With Proposition~\ref{prop:intersec} we
give sufficient conditions to ensure intersection of limit cycles with
one or both lines $x=x_1$ and $x=x_2$. Our result will then follow joining
these informations with the result of Theorem $1$~\cite{Villari85}.

First of all we stress that assumptions are quite standard ones. Hypotheses (A), (B), (C) and (E) guarantee
existence of limit cycles as it will be shown in \S~\ref{ssec:new21}. Hypotheses on $F$
and the equality for $G$ at roots of $F(x)=0$ are {\em fundamental} for our proof.
While we can already find in literature such hypotheses of $F$, the link
between zeros of $F$ and values of $G$ at these points are new, as far as we know.
We remark that hypothesis (C) can be weakened by allowing $F$ to have zeros inside $(x_2,x_1)$, other
than $x_0=0$, where it doesn't change sign.

We already gave some bibliography of results concerning existence
and/or uniqueness of limit cycles for Li\'enard equations; we don't
try to compare our result with all the existing 
ones, we will restrict ourselves to emphasize the strong point of our Theorem
and to compare it with some general results.

First of all we {\em don't assume any parity conditions} on $F$ and/or $g$,
on the contrary if $F$ and $g$ are odd, then Theorem~\ref{thm:main} contains
the Levinson--Smith result~\cite{LevinsonSmith} as particular case: let 
$x_1=-x_2$ be the non--zeros root of $F(x)=0$, $G(x)$ is even because of
oddness of $g$, and then $G(x_1)=G(-x_2)$.

The monotonicity on $F$ is required only outside the interval
determined by the smallest and largest zeros, namely its derivative
$F^{\prime}(x)=f(x)$ can have several zeros inside this interval, this
is a more general situation than the results of Massera~\cite{Massera} and 
Sansone~\cite{Sansone}. The last one follows from our result by
remarking that if $g(x)=x$, then $G(x)=x^2/2$ and let $\Delta >0$
be such that $F(\Delta)=F(-\Delta)=0$, we get $G(\Delta)=G(-\Delta)$.

The Second remark concerns the hypothesis (D): it's easy to verify
if this condition on $G$ holds, just compare the function at two points.
 We don't need to use
the inversion of any function as in the Filippov case~\cite{Filippov}
(and in all results inspired by his method), or to
impose conditions on functions obtained by composition and inversion. These facts 
make our Theorem easily applicable as results of section~\ref{sec:appl} will show.

\section{Main result}
\label{sec:main}

The aim of this section is to prove our main result, Theorem~\ref{thm:main}. 
The proof is divided in two steps, presented in \S~\ref{ssec:parti}
and \S~\ref{ssec:partii}. Before let us introduce two preliminary results, 
first, Proposition~\ref{prop:intersec}, whose role is to give information about
the geometry of limit cycles w.r.t. lines $x=x_i$, where $x_i$ are non zero roots
of $F(x)=0$. Second, give a proof (\S~\ref{ssec:new21}) of existence of limit cycles 
assuming hypotheses (A), (B), (C) and (E), as claimed in the introduction.

\begin{proposition}
\label{prop:intersec}
Let $f$ and $g$ verify hypotheses (A) and (B). Let $F(x)=\int_0^x f(\xi) \, d\xi$
, $G(x)=\int_0^x g(\xi) \, d\xi$ and assume $F(x)$ verify hypothesis (C). Then
\begin{itemize}
\item if $G(x_1) \geq G(x_2)$ all (eventual) limit cycles of~\eqref{eq:lienard2}
will intersect the line $x=x_2$;
\item whereas if $G(x_1) \leq G(x_2)$ all (eventual) limit cycles of~\eqref{eq:lienard2}
will intersect the line $x=x_1$.
\end{itemize}
\end{proposition}
\noindent{\it Proof.\quad}Let us denote by $X_{\mathcal{L}}(x,y)=(y-F(x),-g(x))$ the Li\'enard field 
associate to~\eqref{eq:lienard2} and
let us consider the family of ovals given by: $\mathcal{E}_N=
\{ (x,y)\in \mathbb{R}^2:y^2/2+G(x)-N=0 \}$. 

Let us consider the case $G(x_1) \geq G(x_2)$, the other can be handle 
similarly and we will omit it. The oval $\mathcal{E}_{G(x_2)}$ doesn't
intersect the line $x=x_1$, whereas $\mathcal{E}_{G(x_1)}$ passes
through points $\left(x_2,\pm\sqrt{2\left( G(x_1)-G(x_2)\right)}\right)$. Namely
 $\mathcal{E}_{G(x_1)}$ contains in its interior $\mathcal{E}_{G(x_2)}$
which contains the origin in its interior.

The flow of Li\'enard system~\eqref{eq:lienard2} is transversal to 
$\mathcal{E}_{G(x_2)}$ (more precisely
 it points outward w.r.t to $\mathcal{E}_{G(x_2)}$):
\begin{equation*}
<\nabla \mathcal{E}_{G(x_2)},X_{\mathcal{L}}(x,y)\Big \rvert_{\mathcal{E}_{G(x_2)}}>=-F(x)g(x) \geq 0 \, ,
\end{equation*}
equality holds only for $x=0$ and $x=x_2$. Let us call $(x_1^*,0)$ the unique 
intersection point of $\mathcal{E}_{G(x_2)}$ with the positive $x$--axis.

Hence from Poincar\'e--Bendixson Theorem no limit 
cycle can be completely contained in the strip $[x_2,x^*_1)\times \mathbb{R}$,
moreover orbits of~\eqref{eq:lienard2} spiral outward leaving $\mathcal{E}_{G(x_2)}$.
Thus any (eventual) limit cycle must intersect the line $x=x_2$.\hfill $\Box$

\subsection{Existence of limit cycles}
\label{ssec:new21}

Let us investigate the existence of limit cycles. Consider assumption (E), if 
$\lim_{x\rightarrow \pm \infty}G(x)=+\infty$, we observe that assumption (C) guarantees
that exists $\epsilon >0$ and $\alpha<0<\beta$ such that 
$\int_{\alpha}^{\beta} f(\xi) \, d\xi > \epsilon$. Moreover $f(x) >0$ for $x\not\in[\alpha,\beta]$.
We can then apply Theorem 1 of~\cite{Villari83} to obtain existence of limit cycles.

On the other hand, let us assume $\lim_{x\rightarrow +\infty}G(x)<+\infty$ (the case 
$\lim_{x\rightarrow -\infty}G(x)<+\infty$ can be handle similarly and we omit it), then using
Theorem 3 of~\cite{Villari87a} we complete the proof of the existence of limit cycles.
 
\subsection{Uniqueness: Step I}
\label{ssec:parti}
In~\cite{Villari85} the following result has been proved

\begin{theorem}
\label{thm:villari}
  Let $f$ and $g$ verify hypotheses (A), (B) and  let $F$ verify hypothesis (C). Let $x_2<0<x_1$ be the non--zero roots
of $F(x)=0$. Assume that all limit
cycles of~\eqref{eq:lienard2} intersect the lines $x=x_2$ and $x=x_1$. Then system~\eqref{eq:lienard2} has at most
one limit cycle, if it exists it is stable.
\end{theorem}

Let us give by completeness its proof.

\noindent{\it Proof.\quad}We claim that for any limit cycle, $\gamma$, of system~\eqref{eq:lienard2} we have:
\begin{equation*}
  \oint_{\gamma} g(x) \, dt =0, \quad \oint_{\gamma} g(x)y \, dt =0 \text{ and } 
\oint_{\gamma} g(x)\left[y-F(x)\right] \, dt =0 \, ;
\end{equation*}
this can be proved easily by remarking that $g(x)y=\frac{d}{dt} \left(\frac{1}{2} y^2\right)$. Hence:
\begin{equation}
\label{eq:thmvill1}
\oint_{\gamma} g(x)F(x) \, dt =0 \, .
\end{equation}
Hypotheses (B) and (C) give $F(x)g(x)<0$ for all $x\in (x_2,0)\cup(0,x_1)$, then using the monotonicity of 
$F$ outside $[x_2,x_1]$ and the 
hypothesis that all limit cycles intersect both line $x=x_1$ and $x=x_2$, 
we claim that if $\gamma_1$ and $\gamma_2$ are two limit cycles of~\eqref{eq:lienard2}, $\gamma_1$ contained in the
interior of $\gamma_2$, one has:
\begin{equation*}
  \oint_{\gamma_1} g(x)F(x)\, dt < \oint_{\gamma_2} g(x)F(x)\, dt \, ,
\end{equation*}
which contradicts~\eqref{eq:thmvill1} and so the number of limit cycles is at 
most one.\hfill $\Box$

The weak point of this result is the assumption that all limit cycles must intersect both lines $x=x_1$ and $x=x_2$, in general this 
is not true and moreover it can be difficult to verify. With our result we give sufficient hypotheses to ensure 
this fact. Our Theorem is based
on a slightly generalization of Theorem~\ref{thm:villari} that we state here without proof, which can be obtained 
following closely the previous one.
\begin{theorem}
\label{thm:villarigen}
Assume (A), (B) and (C) of Theorem~\ref{thm:main} hold,
let $N_{x_1,x_2}$ denote the number of limit cycles of system~\eqref{eq:lienard2} 
which intersect both lines $x=x_i$, $i=1,2$. Then $N_{x_1,x_2}\leq 1$.
\end{theorem}
We are now able to prove the main part of our result.

\subsection{Uniqueness: Step II}
\label{ssec:partii}

The number of limit cycles of system~\eqref{eq:lienard2} is by definition
$N_{l.c.}=N_{x_1,x_2}+N_{x_1}+N_{x_2}$, being $N_{x_i}$ the number of limit cycles which
intersect only the line $x=x_i$. So to prove our main result we only need to
control $N_{x_i}$.

>From Proposition~\ref{prop:intersec} and assumption (D) we know that all limit
cycles must intersect {\em both} lines $x=x_i$, $i=1,2$. Namely $N_{x_i}=0$, $i=1,2$.

As already remarked in \S~\ref{ssec:new21} our hypotheses imply existence of at least one limit 
cycle, $N_{l.c.}\geq 1$,
thus we finish our proof recalling that Theorem~\ref{thm:villarigen} gives
$N_{l.c.}\leq 1$.

Before passing to the applications of our Theorem, let us consider in the next
paragraph what can happen when we do not assume hypothesis (D). 

\subsection{Removing the assumption $G(x_1)= G(x_2)$}
\label{ssec:remove4}

The first remark is that assumption (D) cannot be removed
without avoiding cases with more than one limit cycle, as the following
example shows. 
\begin{remark}[A case with $G(x_1)< G(x_2)$]
  Starting from a classical counterexample of Duff and
  Levinson~\cite{LevinsonDuff} to the H. Serbin conjecture~\cite{Serbin},
we exhibit a  polynomial system where all hypotheses (A)--(E) are verified but
  (D), which has $3$ limit cycles. 

Let us consider the equation:
  \begin{equation}
    \ddot x +\epsilon f(x) \dot x + g(x)=0 \, ,
    \label{eq:rem1}
  \end{equation}
where $\epsilon$ is a small parameter, $g(x)=x$ and $f$ is a polynomial of
degree $6$, $f(x)=\sum_{l=0}^{3} a_{2l}x ^{2l}+Ax+Bx^3$, where
$a_0I_0=-4/81$, $a_2I_2=49/81$, $a_4I_4=-14/9$, $a_6I_6=1$,
$I_{2k}=\int_0^{2\pi}\sin^2 \theta cos^{2k}\theta \, d \theta$ and
$A,B$ to be determined. Coefficients $(a_{2l})_{2l}$ are fixed in such a
way that, passing to polar coordinates, for $\epsilon$ small enough and
$A,B=0$, system~\eqref{eq:rem1} has three limit cycles. 

In fact let us introduce polar coordinates: $x=r\cos \theta$, 
$y=r\sin \theta$, then~\eqref{eq:rem1} can be rewritten as:
\begin{equation*}
\begin{cases}
\dot x = y \\
\dot y = -g(x)-\epsilon f(x) y \, ,
\end{cases}
\end{equation*}
thus:\begin{equation*}
\frac{dr}{d\theta}=\frac{\epsilon r f(r\cos \theta)\sin^2 \theta}{1+\epsilon f(r\cos \theta) \sin \theta \cos \theta}\, .
\end{equation*}
If $r$ and $|\epsilon|$ are small enough, we can rewrite the previous equation as:
\begin{equation}
\label{eq:rem2}
  \frac{dr}{d\theta}=\epsilon\left[H_0(r,\theta)+\epsilon H_1(r,\theta)+\epsilon^2 H_2(r,\theta,\epsilon) \right]\, ,
\end{equation}
where $H_i$ are analytic functions of $r,\theta$ and $\epsilon$. Let $\rho>0$ and let us denote 
by $r(\theta,\rho,\epsilon)$,
the solution of~\eqref{eq:rem2} with initial datum $r=\rho$, then our system has a limit cycle
if and only if $\rho$ is an isolated positive root of $r(2\pi,\rho,\epsilon)-\rho=0$. Integrating~\eqref{eq:rem2} 
we get:
\begin{equation}
  r(2\pi,\rho,\epsilon)-\rho=\epsilon \bar{F}(\rho)+\epsilon^2 R_2(\rho,\epsilon) \, ,
\label{eq:rem3}
\end{equation}
where $\bar{F}(\rho)=\int_0^{2\pi} \rho f(\rho \cos\theta)\sin^2\theta \, d\theta$ and $R_2(\rho,\epsilon)$ is some 
analytic remainder function. With our choice of $(a_{2l})_{0\leq l \leq 3}$ we obtain:
 $\bar{F}(\rho)=\rho(\rho^2-1/9)(\rho^2-4/9)(\rho^2-1)$, and 
then from~\eqref{eq:rem3} we conclude that if $|\epsilon|$ is sufficiently small, $r(2\pi,\rho,\epsilon)-\rho$ has 
three positive
 isolated simple roots, $\epsilon$--close to $1/3$, $2/3$ and $1$.

The method used to find the number of limit cycle doesn't involve the values of
$A,B$, we claim that we can vary these parameters in such a way $F(x)=\int_0^x
f(\xi) \, d\xi$ verifies hypothesis (C), with $|x_2|>x_1$ and then $G(x)=x^2/2$
doesn't verify hypothesis (D). Just as an example consider:
\begin{equation*}
 F(x)=\frac{x}{\pi}\left(-\frac{4}{81} +\frac{196}{81}\frac{x^2}{3}-\frac{112}{9}\frac{x^4}{5}+
\frac{64}{5}\frac{x^6}{7}+\frac{1}{200} x+\frac{1}{2} x^3\right) \, ,
\end{equation*}
which has three real zeros $x_0=0$, $x_2<0<x_1$ and its monotone increasing outside 
$(x_2,x_1)$. Moreover $f(x)=F^{\prime}(x)$ has four zeros in the same 
interval~\footnote{Using Sturm's method to find real roots of polynomials we obtain that 
the zeros of $F$ belong to the intervals: $x_2\in [-1.130,-1.129]$ and $x_1\in[0.247,0.248]$, whereas  
zeros of $f$ verify: $x_4^{\prime}\in [-0.969,-0.9688]$,$x_3^{\prime}\in[-0.343,-0.342]$, 
$x_2^{\prime}\in [-0.173,-0.172]$ and $x_1^{\prime}\in[0.139,0.140]$.}.
\end{remark}

To conclude this part let us remark that adding further assumptions on $F(x)$,
one can ensure that all limit cycles must intersect both lines 
$x=x_1$, $x=x_2$, thus obtaining a {\em existence and uniqueness} result 
for~\eqref{eq:lienard2}. For instance one can prove the following
\begin{theorem}
\label{thm:furtherass}
Assume $f$ and $g$ verify hypotheses (A) and (B). Let $F$ and $G$ be 
the primitives of $f$ and $g$ vanishing at $x=0$ and assume they verify 
hypotheses (C) and (E). Assume one of the following conditions hold:
 \begin{itemize}
\item[$(D^{\prime})$] $G(x_1) > G(x_2)$ and there exists $x_2^*\in (x_2,0)$
such that $F(x_2^*)\geq \sqrt{2G(x_1)}$;
\item[$(D^{\prime\prime})$] $G(x_1) < G(x_2)$ and there exists $x_1^*\in(0,x_1)$
such that $F(x_1^*)\leq -\sqrt{2G(x_2)}$.
\end{itemize}
Then Li\'enard system~\eqref{eq:lienard2} has one and only one limit cycle. 
\end{theorem}

\noindent{\it Proof.\quad}
We only prove the previous Theorem assuming $(D^{\prime})$, being the other case very similar.
Let us assume $G(x_1) > G(x_2)$ and that there exists $x_2^*\in (x_2,0)$
such that $F(x_2^*)\geq \sqrt{2G(x_1)}$, we will prove that any orbit which intersects the line $x=x_2^*$
must intersect also the line $x=x_1$.

Considering the oval $\mathcal{E}_{G(x_1)}=
\{ (x,y)\in \mathbb{R}^2: y^2/2+G(x)-G(x_1)=0 \}$ one realizes that there exists a unique point
$(0,y_A)$ with $y_A< \sqrt{2G(x_1)}$, whose future orbit will intersect the line $x=x_1$ at the point
$(x_1,0)$.

Let us consider now a point $(x_2^*,y_B)$, with $y_B\geq F(x_2^*)$, we claim that its future orbit
will intersect the $y$--axis at some $(0,y_{B^{\prime}})$ such that $y_{B^{\prime}}> \sqrt{2G(x_1)}$.
This can be proved by considering the evolution of the function $\Lambda(x,y)=y^2/2+G(x)$ under the flow
of the Li\'enard system.

Summarizing the orbit of all point of the form $(x_{2}^{*},y_{B})$ such that 
$y_{B}> F(x_{2}^{*})$, will intersect the line $x=x_1$ with positive $y$ coordinate. This conclude the proof
once we remark that orbits of points $(x_{2}^{*},y^{\prime})$ such that $y^{\prime} < F(x_{2}^{*})$, turn clockwise
and will intersect again the line $x=x_{2}^{*}$ at some point $(x_{2}^{*},y^{\prime\prime})$ with 
$y^{\prime\prime}\geq F(x_{2}^{*})$.

To complete the proof of the Theorem one remark that by Proposition~\ref{prop:intersec} all limit cycles
must intersect the line $x=x_2$, hence they must intersect the line $x=x_2^*$, being $x_2<x_2^*$. By the first part
these limit cycles intersect also the line $x=x_1$ and then applying Theorem~\ref{thm:villari} we conclude 
the proof.\hfill $\Box$

\section{Some applications}
\label{sec:appl}

In this section we give some applications of Theorem~\ref{thm:main}.
The first application concerns Li\'enard's systems~\eqref{eq:lienard2},
where $F$ and $G$ verify all hypotheses of Theorem~\ref{thm:main} but
(D) (\S~\ref{ssec:caseI} and \S~\ref{ssec:caseII}). 
Our aim is to show that we can find a new Li\'enard system 
(slightly modified version of the original one) for which Theorem~\ref{thm:main} holds, then exhibiting one and 
only a limit cycle. The second application is of different nature, starting with a given Li\'enard system,
which doesn't verify assumptions of Theorem~\ref{thm:main}, we prove existence and uniqueness of limit cycles for
a new system obtained from the first one just by introducing two parameters. We will consider the polynomial
case (\S~\ref{ssec:caseP}) and a more general one (\S~\ref{ssec:secT}). 

\subsection{Case I: deform $g$}
\label{ssec:caseI}

Let us recall that $F$ has three real zeros, $x_0=0$ and $x_2<0<x_1$,
let us assume $G(x_1)\neq G(x_2)$. Let us introduce the $1$--parameter
family of functions:
\begin{equation}
\label{eq:glambda}
g_{\lambda}(x)=\begin{cases}
               g(x) & \text{ if $x\geq 0$} \\
               \lambda g(x) & \text{ if $x< 0$} \, .
               \end{cases}
\end{equation}
Then $(g_{\lambda})_{\lambda}$ verifies hypotheses (A) and (B)
of Theorem~\ref{thm:main}, provided $\lambda >0$. Let us define
$G_{\lambda}(x)=\int_0^x g_{\lambda}(\xi) \, d\xi$;
 let $\lambda_*=G(x_1)/G(x_2)>0$, then $G_{\lambda_*}(x_1)= G_{\lambda_*}(x_2)$.
 Hence also hypotheses (D) and (E) hold and the differential equation:
\begin{equation*}
\ddot x + f(x) \dot x +g_{\lambda_*}(x) = 0 \, ,
\end{equation*}
has a unique isolated periodic solution.

\subsection{Case II:  deform $F$}
\label{ssec:caseII}

Let us assume $G(x_1)<G(x_2)$. The idea
is now to modify the roots of $F$ in such a way hypothesis (D) holds.
We do this in a simple way, more sophisticated ones are possibles.

Let $\lambda >0$ and let us introduce the $1$--parameter family of
functions $(F_{\lambda})_{\lambda}$, defined by:
\begin{equation*}
F_{\lambda}(x)=\begin{cases}F(x) & \text{ if $x\geq 0$} \\
F(\lambda x) & \text{ if $x< 0$} \, ,
\end{cases}
\end{equation*}
clearly $(F_{\lambda})_{\lambda}$ verifies hypothesis (E)
if $F$ does; $(F_{\lambda})_{\lambda}$ is no longer Lipschitz at $x=0$
but existence and uniqueness of the Cauchy problem are still verified.

Thanks to the form of $g$ and hypothesis on $G$, there exists a unique
$x_2^{*}<0$ such that $G(x_2^{*})=G(x_1)$, moreover $x_2<x_2^{*}$. Let 
$\lambda_* =\frac{|x_2|}{|x_2^{*}|}$ and $\bar{x}_{\lambda_*}=x_2/\lambda_*$. 
We claim that $\bar{x}_{\lambda_*}$ is the unique negative zeros of
 $F_{\lambda}(x)$. Hence hypotheses (C) and (D) hold, in fact: $F_{\lambda}$ 
has three zeros, $x_0$, $x_1>0$ (as $F$ does) and $\bar{x}_{\lambda}$, moreover
 $G(\bar{x}_{\lambda})=G(x_2^{*})=G(x_1)$. Hence 
  \begin{equation*}
    \begin{cases}
      \dot x =y - F_{\lambda_*}(x) \\
      \dot y = -g(x) \, ,
    \end{cases}
  \end{equation*}
has a unique limit cycle.

\subsection{Polynomial Case}
\label{ssec:caseP}

Let us consider a polynomial
$P_{2n+1}(x)=a_{2n+1}x^{2n+1}+a_{2n}x^{2n}+\dots +a_1 x$, assume
$n\geq 1$, $a_{2n+1}>0$ and hypothesis (C) doesn't hold. We claim that
we can introduce a modified Polynomial
$P_{\lambda}(x)=P_{2n+1}(x)-\lambda x$ and a function $g$ verifying
hypotheses (A), (B) and (D) such that
\begin{equation}
\label{eq:polycase}
  \begin{cases}
    \dot x = y-P_{\lambda}(x) \\
    \dot y = -g(x) \, ,
  \end{cases}
\end{equation}
has a unique limit cycle.

$P_{2n+1}(x)$ has at most $2n$ local maxima and minima, so let us
define:
\begin{eqnarray*}
  \xi_{+} & =\min \{ x>0: \forall y>x : P^{\prime}_{2n+1}(y)>0
\text{ and } P^{\prime\prime}_{2n+1}(y)>0\} \\
\xi_{-} & =\max \{ x<0: \forall y<x : P^{\prime}_{2n+1}(y)>0\text{
  and } P^{\prime\prime}_{2n+1}(y)>0 \} \, .
\end{eqnarray*}
Let us consider $\lambda_{\pm}\geq 0$ such that:
\begin{equation*}
  P_{2n+1}(x)\leq \lambda_{+} x \quad \text{for all $0<x<\xi_{+}$
   and } P_{2n+1}(x)\geq \lambda_{-} x \quad \text{for all $\xi_{-}<x<0$.}
\end{equation*}
Such $\lambda_{\pm}$ can be obtained as follows. Consider straight
lines $y=\mu x$ tangent to $y=P_{2n+1}(x)$ for $x\in (0,\xi_+)$, they
are in finite number, so one can take $\lambda_{+}=\max
|\mu_i|$; if $P_{2n+1}(x)<0$ on $(0,\xi_+)$ we set $\lambda_+=0$. A
similar construction can be done for $\lambda_{-}$.

Let $\bar{\lambda}=\max \{ \lambda_+,\lambda_{-} \}$ then we claim
that for all $\lambda > \bar{\lambda}$,
$P_{\lambda}(x)=P_{2n+1}(x)-\lambda x$ satisfies hypothesis (C). By
construction $P_{\lambda}(x)<0$ for all $x\in (0,\xi_+)$ and
$P_{\lambda}(x)>0$ for all $x\in (\xi_-,0)$. Because $a_{2n+1}>0$, for
sufficiently large $|x|$, $P_{\lambda}(x)$ has the same sign than $x$,
then for $x>0$ large enough, $P_{\lambda}(x)>0$, hence there is at
least one zeros of $P_{\lambda}(x)$. Actually this will be the only
one. For suppose there are more zeros~\footnote{They will be at least
  three, if transversal, because $P_{\lambda}(\xi_+)<0$ and
  $P_{\lambda}(x)>0$ for $x$ large enough. Non--transversal zeros can
  be removed by small increment of $\lambda$.}
and call them $\bar{x}_1<\bar{x}_2<\bar{x}_3$; by construction for all $x\in
(\bar{x}_1,\bar{x}_2)$ we have $P_{2n+1}(x)> \lambda x$ whereas
  $P_{2n+1}(x)< \lambda x$ for $x\in (\bar{x}_2,\bar{x}_3)$. This
    implies $P_{2n+1}(x)$ non--convex for $x>\xi_+$, against the
    definition of $\xi_+$. The case for negative $x$ can be handle in
    a similar way. Let us call $x_1$ the positive zeros and $x_2$ the
    negative one. 
Summarizing:
$P_{\lambda}(x)$ has threes real zeros: $x_0=0$, $x_2<0<x_1$, moreover
$P_{\lambda}(x)<0$ for $0<x<x_1$, and $P_{\lambda}(x)>0$ for
$x_2<x<0$. Remark that $x_1> \xi_+$ and $x_2<\xi_-$, namely
$P_{\lambda}(x)$ is monotone increasing outside $[x_2,x_1]$.

Let $g$ be any locally Lipschitz function such that: $xg(x)>0$ for $x\neq
0$ and $\int_{x_2}^{x_1} g(\xi) \, d\xi = 0$, then
Theorem~\ref{thm:main} applies and~\eqref{eq:polycase} has a unique
limit cycle.

\subsection{Generalization of the polynomial case.}
\label{ssec:secT}

In this section we will generalize the result of the previous section,
by proving an existence and uniqueness result for
the Li\'enard equation.

\begin{theorem}
 Let us consider the Li\'enard equation:
\begin{equation}
 \ddot x +f(x) \dot x +g(x)=0 \, ,
 \label{eq:them2}
\end{equation}
where $f$ and $g$ verify:
\begin{itemize}
  \item[(A)] $f$ is continuous and $g$ is locally Lipschitz;
  \item[($B^{\prime}$)] $\lim_{x\rightarrow \pm \infty} f(x)=+\infty$ and $xg(x)>0$ for all $x\neq 0$.
\end{itemize}
Then there exist $\hat{\lambda}$, such that for
all $\lambda \geq \hat{\lambda}$ there exists $\mu=\mu(\lambda)$ and system:
  \begin{equation}
    \ddot x +f_{\lambda}(x) \dot x +g_{\mu}(x) =0\, ,
    \label{eq:them22}
  \end{equation}
has a unique limit cycle, where $f_{\lambda}(x)=f(x)-\lambda$ and
$g_{\mu}$ will be defined in~\eqref{eq:defgmu}.
\end{theorem}

\begin{remark}
  Hypothesis ($B^{\prime}$) is a strong one, even if it is always
  verified for the important class of polynomial Li\'enard equations. It
  can be relaxed by assuming $\lim_{x\rightarrow \pm \infty}F(x)=\pm
  \infty$ and $F$ to be monotone increasing outside some interval containing the origin, 
 where as usual $F(x)=\int_0^x f(\xi)
  \, d\xi$. 
\end{remark}

\noindent{\it Proof.\quad}
For any $\lambda_1>f(0)$, system~\eqref{eq:them2} where
$f_{\lambda_1}(x)=f(x)-\lambda_1$ replaces $f(x)$, has at least a
limit cycle (see Theorem 3 of~\cite{Villari87a}). Then one can find
a $\hat{\lambda}\geq \lambda_1$ such that for all $\lambda \geq
\hat{\lambda}$,
$F_{\lambda}(x)=-\lambda x+ \int_0^x f(\xi) \, d\xi$ verifies
hypotheses of Theorem~\ref{thm:main}. Just use
monotonicity of $F$, as we did in the previous section for the polynomial
case, to ensure that with $\lambda$ large enough, $F_{\lambda}$ has only two non zeros roots and
it is monotone increasing outside the interval whose boundary is formed by the two non zeros roots.. 

Let us call $x_2(\lambda)<0<x_1(\lambda)$, the non--zeros roots of
$F_{\lambda}(x)=0$. Then we can modify $g$ (for instance as we did
in~\S~\ref{ssec:caseI}) introducing:
\begin{equation}
\label{eq:defgmu}
g_{\mu}(x)= \begin{cases}
            g(x) &\text{ if $x\geq 0$}; \\
            \mu g(x) &\text{ if $x< 0$},
  \end{cases}
\end{equation}
in such a way $\int_{x_2}^{x_1} g(\xi) \, d\xi = 0$. Namely also
hypothesis (D) of Theorem~\ref{thm:main} holds, and so
system~\ref{eq:them22} has a unique limit cycle.\hfill $\Box$

The role of $f$ and $g$ in the previous Theorem can be in some sense inverted.
More precisely, one can prove the following result
\begin{remark}
Let us given the global center system:
\begin{equation}
\label{eq:cent}
\ddot x + g(x)=0 \, ,
\end{equation}
with $g$ locally Lipschitz, $xg(x) >0$ for $x\neq 0$, 
$G(x)= \int_0^x g(\xi) \, d\xi$ and assume  
$\lim_{x\rightarrow \pm \infty}G(x) = +\infty$. 
Take any $x_2<0<x_1$ such that $G(x_2)=G(x_1)$ then we can 
{\em perturb}~\eqref{eq:cent} by adding {\em any continuous friction term} 
$f(x)\dot x$, such that $F(x)=\int_0^x f(\xi) \, d\xi$ verifies 
$F(x_1)=F(x_2)=0$ and $F(x)$ is monotone increasing outside the interval $[x_2,x_1]$,
obtaining a Li\'enard system: $\ddot x +f(x)\dot x+g(x)=0$ with one and only
one limit cycle.
\end{remark}

\end{document}